\documentclass[11pt]{amsart}

\usepackage{amsmath,amsthm, amscd, amssymb, amsfonts}
\usepackage{epsfig}
\input xy
\xyoption{all}

\newcommand{\Mo}{{\mathcal M}}
\newcommand{\No}{{\mathcal N}}

\newcommand{\Ss}{{\mathcal S}}
\newcommand{\ot}{{\otimes}}
\newcommand{\otk}{{\ot_{\ku}}}
\newcommand{\otb}{{\overline{\otimes}}}
\newcommand{\otc}{\widetilde{{\otimes}}}

\newcommand{\ele}{{\mathcal L}}

\newcommand{\ca}{{\mathcal C}}

\newcommand{\ac}{{\mathcal A}}
\newcommand{\Fc}{{\mathcal F}}
\newcommand{\Gc}{{\mathcal G}}

\newcommand{\cop}{\rm{cop}}

\newcommand{\ku}{{\Bbbk}}
\newcommand{\ax}{{\mathbb A}}

\newcommand{\uno}{{\mathbf 1}}
\newcommand{\C}{{\mathbb C}}

\newcommand{\id}{\mbox{\rm id\,}}

\newcommand{\Ind}{\mbox{\rm Ind\,}}

\newcommand\Rep{\operatorname{Rep}}

\newcommand\co{\operatorname{co}}

\newcommand\Hom{\operatorname{Hom}}

\newcommand{\End}{\operatorname{End}}

\newcommand\St{\operatorname{Stab}}

\newcommand{\Id}{\mathop{\rm Id}}

\renewcommand{\_}[1]{\mbox{$_{\left( #1 \right)}$}}

\theoremstyle{plain}

\numberwithin{equation}{section}

\newtheorem{teo}{Theorem}[section]

\newtheorem{lema}[teo]{Lemma}

\newtheorem{prop}[teo]{Proposition}

\theoremstyle{definition}

\newtheorem{defi}[teo]{Definition}

\theoremstyle{remark}

\newtheorem{rmk}[teo]{Remark}

\def\pf{\begin{proof}}

\def\epf{\end{proof}}

\theoremstyle{remark}

\begin{document}

\title[Constructing dynamical twists ]
{Constructing dynamical twists over a non-abelian base}
\author[ Mombelli]{
Juan Mart\'\i n Mombelli }
\thanks{
\newline
\emph{keywords}: Dynamical Yang-Baxter equation; dynamical twists;
Hopf algebras}
\address{Mathematisches Institut, Universit\"{a}t M\"{u}nchen, Theresienstra\ss e 39,
\newline \indent D-80333 M\"{u}nchen, Germany.} \email{mombelli@mate.uncor.edu
\newline \indent\emph{URL:}\/ http://www.mate.uncor.edu/mombelli}
\begin{abstract} We give examples of dynamical twists in
finite-dimensional Hopf algebras over an  arbitrary Hopf
subalgebra. The construction is based on the categorical approach
of dynamical twists introduced by Donin and Mudrov \cite{DM1}.

\end{abstract}

\date{\today}
\maketitle

\setcounter{tocdepth}{1}
\section*{Introduction}

The theory of \emph{dynamical quantum groups} initiated by G.
Felder \cite{F1} is nowadays an active branch of mathematics. This
theory arose from the notion of \emph{dynamical Yang-Baxter
equation}, also known as the \emph{Gervais-Neveu-Felder equation}
in connection with integrable models of conformal field theories
and Liouville theory, see for example \cite{F2}, \cite{GN},
\cite{ABB}. For a detailed review and bibliography on the
dynamical Yang-Baxter equation the reader is referred to \cite{E},
\cite{ES}, \cite{Fh}.

\medbreak

In \cite{B}, \cite{BBB}, the notion of Drinfeld's twist for Hopf
algebras was generalized in the dynamical setting. A
\emph{dynamical twist} for a Hopf algebra $H$ over $A$, where $A$
is an abelian subgroup of the group of grouplike elements in $H$,
is a function $J(\lambda):\widehat{A}\to (H\ot H)^{\times}$
satisfying certain non-linear equations. When the group $A$ is
trivial we recover the notion of Drinfeld's twist. When $H$ is a
quasi-triangular Hopf algebra with $R$-matrix $R$ then
$\mathcal{R}(\lambda)=J^{-1}(\lambda)^{21} RJ(\lambda)$ satisfies
the dynamical quantum Yang-Baxter equation. See \cite{EN1}.

\medbreak

In the finite-dimensional case, dynamical twists were studied
first in \cite{EN1} and later in \cite{M}. In the first paper
dynamical twists over an abelian group for the group algebra of a
finite group are classified. Following closely this work, in
\cite{M}, dynamical twists over an abelian group for any
finite-dimensional Hopf algebra are described.

\medbreak

In this work we extend the results of \cite{M} in the case where
$A$ is an arbitrary Hopf subalgebra of $H$. To this end we rely on
the definition and categorical construction of dynamical twists
introduced in \cite{DM1} from dynamical adjoint functors. Here,
the Hopf subalgebra $A$ plays the same role as the Levi subalgebra
of a reductive Lie algebra in \cite{DM1}. As in \cite{M}, the
language of \emph{module categories} has been used with profit.

\medbreak

The contents of the paper are organized as follows: in Section 1
we recall the notion of stabilizers for Hopf algebra actions
introduced by Yan and Zhu \cite{YZ}, and the definition of module
categories over a tensor category.

\medbreak

In Section 2 we give a brief account of the results and
definitions appeared in \cite{DM1}. We explain the definition of
dynamical extension of a tensor category,  dynamical twist and
dynamical adjoint functors. We also recall the construction of a
dynamical twist coming from a pair of dynamical functors.

\medbreak

Section 3 is devoted to the construction of dynamical twists for a
finite-dimensional Hopf algebra $H$ over a Hopf subalgebra $A$.
Following \cite{EN1} we shall say that a \emph{dynamical datum }
for $(H, A)$ is a pair $(K,T)$, where
\begin{itemize}
    \item $K$ is a $H$-simple left $H$-comodule algebra with $K^{\co H}=\ku$,
    \item and $T:\Rep(A)\to {}_K\Mo$ is a functor such that  \end{itemize}
for any $V, W\in \Rep(A)$ there are natural isomorphisms
$$\St_K(T(V),T(W))\simeq \big(\Ind^H_A\, (V\ot_{\ku}
W^*)\big)^*.$$ For any dynamical datum $(K,T)$ we shall show that
the pair $(T, R)$, where $R:\Rep(H)\to  \Rep(A)$ is the
restriction functor, is a pair of dynamical adjoint functors.
Therefore, applying the tools explained in Section 2, we obtain a
dynamical twist over the base $A$. This construction generalize
the procedure appeared in \cite{M} when $A$ is a commutative
cocommutative Hopf subalgebra. Also, we shall prove that any
dynamical twist based in $A$ for $H$ comes from a dynamical datum.
 Finally we show some explicit examples of
dynamical data.

\section{Preliminaries and notation}

Throughout this paper $\ku$ will denote an arbitrary field. All
categories and functors are assumed to be $\ku$-linear. All vector
spaces and algebras are assume to be over $\ku$. If $K$ is an
algebra, we shall denote by ${}_K\Mo$ the category of
finite-dimensional left $K$-modules.

If $V$ is a vector space, we shall denote by $\langle\;,\,
\rangle:V^*\ot_{\ku} V\to \ku$ the evaluation map.

By $H$ we shall denote a Hopf algebra with counit $\varepsilon$,
and antipode $\Ss$. We shall denote by $\Rep(H)$ the category of
finite-dimensional left $H$-modules emphasizing the canonical
tensor structure.

If $K$ is an $H$-comodule algebra with coaction $\delta:K\to H\ot
K$, an $H$-\emph{costable ideal} of $K$ is a two-sided ideal $I$
of $K$ such that $\delta(I)\subseteq H\ot I$. We shall say that
$K$ is \emph{$H$-simple} if it has no non-trivial $H$-costable
ideal of $K$.

We shall denote $^{H}\Mo_K$ the category of left $H$-comodules,
right $K$-modules such that the $K$-module structure is an
$H$-comodule map. If $P \in {}^{H\!}\Mo_K$ then $\End_K(P)$ has a
natural left $H$-comodule algebra structure via
$\delta:\End_K(P)\to H\ot_{\ku} \End_K(P) $, $T\mapsto T\_{-1}\ot
T\_0$, determined by
\begin{equation}\label{h-comod} \langle\alpha, T\_{-1}\rangle\,
T_0(p)=\langle\alpha, T(p\_0)\_{-1}\Ss^{-1}(p\_{-1})\rangle\,
T(p\_0)\_0,\end{equation}  $T\in\End_K(P),$ $p\in P$, $\alpha\in
H^*$. See \cite[Lemma 1.26]{AM}.

\begin{lema}\label{preliminar1} Let $A\subseteq H$ be a Hopf
subalgebra and $V$ an $A$-module. The space $\Hom_A(H, V)$ has a
natural $H$-module structure and  there are natural $H$-module
isomorphisms
$$\big(\Ind^H_A\, V \big)^*\simeq \Hom_A(H,  V^*).$$
\end{lema}

\pf The $H$-module structure on $\Hom_A(H,  V)$ is as follows. If
$t, h\in H$, $T\in \Hom_A(H,  V)$ then
$$(h\cdot T)(t)=T(th).$$

It is not difficult to prove that the maps
$$\theta:\big(\Ind^H_A\, V \big)^*\to \Hom_A(H,
V^*),$$
$$\widetilde{\theta}: \Hom_A(H,
V^*)\to \big(\Ind^H_A\, V\big)^*,$$ given by the formulas
\begin{align*} \theta(\alpha)(h)&=\sum_{i} \alpha(\overline{\Ss(h)\ot v_i})\,\, v^i,\\
\widetilde{\theta}(\beta)(\overline{h\ot
v})&=\langle\beta(\Ss^{-1}(h)), v\rangle,
\end{align*}
are well defined isomorphisms, one the inverse of each other, and
they are $H$-module maps. Here $ \alpha\in \big(\Ind^H_A\, V
\big)^*, \beta\in\Hom_A(H,  V^*)$, and $(v_i), ( v^i)$ are dual
basis for $V$ and $V^*$, $h\in H,$ $\overline{h\ot v}\in
\Ind^H_A\, V$.

\epf

\subsection{Stabilizers for Hopf algebra actions}  We recall very
briefly the notion of Hopf algebra stabilizers introduced in
\cite{YZ}, see also \cite{AM}.

\medbreak

Let $K$ be a finite-dimensional left $H$-comodule algebra and $V,$
$W$ two left $K$-modules. The \emph{Yan-Zhu stabilizer}
$\St_K(V,W)$ is defined as the intersection
$$ \St_K(V,W)= \Hom_K(H^* \ot V,H^* \ot W)
\cap \ele\big(H^* \ot \Hom (V, W)\big).$$ Here the map $\ele:H^*
\ot \Hom (V, W)$ $\to$ $ \Hom(H^* \ot V,H^* \ot W)$ is defined by
$\ele(\gamma\ot T)(\xi\ot v)=\gamma\xi\ot T(v)$, for every
$\gamma, \xi\in H^*$, $T\in \Hom (V, W)$, $v\in V$.

\medbreak

The $K$-action on $H^*\ot V$ is given by
$$ k\cdot (\gamma\ot v)= k_{(-1)} \rightharpoondown \gamma\, \ot\,
k_{(0)}\cdot v,$$ for all $k\in K$, $\gamma\in H^*$, $v\in V$.
Here $\rightharpoondown:H\ot H^*\to H^*$ is the action defined by
$\langle h\rightharpoondown \gamma, t\rangle=\langle \gamma,
\Ss^{-1}(h)t\rangle$, for all $h, t\in H$, $\gamma\in H^*$. Also,
we denote $\St_K(V)=\St_K(V,V)$.

\begin{prop}\cite[Prop. 2.7, Prop. 2.16]{AM}\label{stab-properties} The following assertions holds.

1. For any left $K$-modules $V, W, U$ there is a natural
composition
$$ \St_K(V,W)\ot_{\ku} \St_K(U, V)\to \St_K(U,W)$$
making $\St_K(V)$  a left $H$-module algebra.

2. If $K$ is $H$-simple then
\begin{equation}\label{dim-stab} \dim(K) \dim(\St_K(V, W))= \dim(V)
\dim(W) \dim(H).
\end{equation}

3. For any $X\in \Rep(H)$ there are natural isomorphisms
$$\Hom_H(X,\St_K(V, W))\simeq \Hom_K(X\ot_{\ku} V, W),$$
where the action on $X\ot_{\ku} V$ is given by the coaction of
$K$.\qed
\end{prop}

The following result concernig Yan-Zhu stabilizers will be useful
later. If $A\subseteq H$ is a Hopf subalgebra, and $R=K^{\co
A}\subseteq K$ is a left $A$-Hopf Galois extension then there are
$H$-module algebra isomorphisms
\begin{align}\label{stab-galois} \St_K(V,W)\simeq
\Hom_{A}(H, \Hom_R(V, W))
\end{align}
for any left $K$-modules $V, W$. Worth to mention that the action
of $A$ on $\Hom_R(V, W)$ is given by
$$ (a\cdot T)(v)= a^{[1]}\cdot T(a^{[2]}\cdot v),$$
for all $a\in A, T\in\Hom_R(V, W), v\in V$. Recall that the map
$\gamma:A\to K\ot_R K$, $\gamma(a)=a^{[1]}\ot a^{[2]}$ is defined
by $\gamma(a)= can^{-1}(a\ot 1)$, where $can:K\ot_R K\to
A\ot_{\ku} K$ is the canonical map $can(k\ot s)=k\_{-1}\ot k\_0
s$, $k,s \in K$. For more details see \cite[Rmk. 3.4]{Sch},
\cite[Thm. 2.23]{AM}.

\subsection{Module categories}\label{subsection-modcat} We briefly recall the definition of module
category and the definition introduced by Etingof-Ostrik of
\emph{exact} module categories. We refer to \cite{O1}, \cite{O2},
\cite{eo}.

\medbreak

Let us fix $\ca$ a finite tensor category. A \emph{module
category} over $\ca$ is a collection $(\Mo,\otb, m, l)$ where
$\Mo$ is an Abelian category, $\otb:\ca\times\Mo\to \Mo$ is an
exact bifunctor, associativity and unit isomorphisms
$m_{X,Y,M}:(X\ot Y)\otb M\to X\otb ( Y\otb M)$, $l_M:\uno\otb M\to
M$, $X,Y\in \ca$, $M\in\Mo$, such that
\begin{align}\label{modcat1} m_{X,Y,Z\otb M}\, m_{X\ot
Y,Z,M}&=(\id_X\ot\, m_{Y,Z,M})\,m_{X,Y\ot Z,M}\,(a_{X,Y,Z}\ot\id_M),\\
\label{modcat2} (\id_X\ot\, l_M)\, m_{X,\uno,M}&= r_X\ot\id_M,
\end{align}
for all $X,Y, Z\in \ca$, $M\in\Mo$. Sometimes we shall simply say
that $\Mo$ is a module category omitting to mention $\otb, m $ and
$l$ whenever no confusion arises.

\medbreak

In this paper we further assume that all module categories have
finitely many isomorphism classes of simple objects.

\medbreak

Let $\Mo$, $\Mo'$ be two module categories over $\ca$. A module
functor between $\Mo$ and $\Mo'$ is a pair $(\Fc,c)$ where
$\Fc:\Mo\to \Mo'$ is a functor and $c_{X,M}:\Fc(X\otb M)\to X\otb
\Fc(M)$ is a family of natural isomorphisms such that
\begin{align}\label{modfun1} m'_{X,Y,\Fc(M)}\, c_{X\ot Y,M}&=
(\id_X\ot\,
c_{Y,M})\, c_{X,Y\otb M} \Fc(m_{X,Y,M}),\\
\label{modfun2} l'_{\Fc(M)}\, c_{\uno, M} &= \Fc(l_M),
\end{align}
for all $X,Y\in \ca$, $M\in\Mo$. If $(\Gc,d):\Mo\to \Mo'$ is
another module functor, a morphism between $(\Fc,c)$ and $(\Gc,d)$
is a natural transformation $\alpha: \Fc \to \Gc$ such that for
any $X\in \ca$, $M\in \Mo_1$:
\begin{gather}
\label{modfunctor3} d_{X,M}\alpha_{X\otimes M} =
(\id_{X}\ot\alpha_{M})c_{X,M}.
\end{gather}

The module structure over $\Mo\oplus \Mo'$ is defined in an
obvious way. A module category is \emph{indecomposable} if it is
not equivalent to the direct sum of two non-trivial module
categories.

A module category $\Mo$ is \emph{exact} \cite{eo} if for any
projective object $P\in \ca$ and any $M\in \Mo$ the object $P\otb
M$ is again projective.

\subsection{Module categories over Hopf algebras}\label{cathopf} Let $H$ be a
finite-dimensional Hopf algebra. Let $K$ be a left $H$-comodule
algebra. Then ${}_K\Mo$ is a left module category over $\Rep H$
via the coaction $\lambda: K \to H\otimes K$. That is, $\ot:\Rep
H\times{}_K\Mo\to {}_K\Mo$ is given by
$$X\ot V:=X\ot_{\ku} V,$$ for $X\in\Rep H$ and $V\in{}_K\Mo$
with  action $k\cdot (x\ot v)=k\_{-1}\cdot x\,\ot \, k\_0\cdot v$,
for all $k\in K$, $x\in X$, $v\in V$. Moreover, any exact module
category is of this form.

\begin{teo}\label{cathopf-resultados} \begin{enumerate}
    \item[1.] If $K$ is $H$-simple left $H$-comodule algebra then
    ${}_K\Mo$ is an indecomposable exact module category.
    \item[2.] If $\Mo$ is an indecomposable exact module category over
    $\Rep(H)$ then there exists an $H$-simple left $H$-comodule
    algebra $K$ such that $\Mo\simeq {}_K\Mo$.
\end{enumerate}
\end{teo}
\pf See \cite[Prop. 1.20, Th. 3.3]{AM}.\epf

Let $S$ be another $H$-simple left $H$-comodule algebra.

\begin{prop}\label{eq1}\cite[Prop. 1.24]{AM} The module categories ${}_K\Mo$, ${}_S\Mo$
over $\Rep(H)$ are equivalent if and only if there exists $P\in
{}^{H\!}\Mo_{K}$ such that $S\simeq \End_K(P_K)$ as $H$-module
algebras. Moreover the equivalence is given by $F:{}_K\Mo\to
{}_S\Mo$, $F(V)=P\ot_K V$, for all $V\in {}_K\Mo$. \qed
\end{prop}

\section{Dynamical twists constructed from dynamical functors}
We recall a construction due to Donin and Mudrov of dynamical
twists from dynamical adjoint functors, see \cite[\S 6]{DM1}.
There are some differences in our statements and those appeared in
\emph{loc. cit.} since we use left module categories instead of
right ones.

\subsection{Dynamical extensions of tensor categories}
In \cite{DM1} for any  tensor category $\ca$ and a module category
$\Mo$ over $\ca$ the authors introduced a new tensor category,
that we will denote by $\Mo\ltimes \ca$. This tensor category is
called the \emph{dynamical extension} of $\ca$ over $\Mo$.

\medbreak

Objects in the category $\Mo\ltimes \ca$ are functors $F_X:\Mo\to
\Mo$, $F_X(M)=X\overline{\ot} M$, for all $X \in \ca$, $M\in \Mo$.
Morphisms are natural transformations. Observe that for each $f\in
\Hom_{\ca}(X, Y)$ there is a natural transformation $\eta_f:F_X\to
F_Y$, given by
$$(\eta_f)_M: X\otb M\to Y\otb M,\quad  (\eta_f)_M=
f\ot\,\id_M,$$ for all $M\in \Mo$.

We briefly recall the monoidal structure of $\Mo\ltimes \ca$. The
tensor product is $F_X\ot F_Y=F_{X\ot Y}$, $X, Y\in \ca$, and the
associativity constraint is
$$\widetilde{a}_{X,Y,Z}:(F_X\ot F_Y)\ot F_Z\to  F_X\ot (F_Y\ot
F_Z), \quad  (\widetilde{a}_{X,Y,Z})_M=(a_{X,Y,Z}\ot \id_M),$$ for
all $M\in \Mo$. For any $X\in \ca$ the left and right unit
isomorphisms are given by
$$ l_X:F_X\ot F_{\uno}\to F_X,\quad r_X:F_{\uno}\ot F_X\to F_X,$$
where $l_{X,M}=l_X\ot\id_M$ and $r_X=r_X\ot\id_M$ for all $M\in
\Mo$.

\medbreak

If $\eta:F_X\to F_Z$, $\phi:F_Y\to F_W$ are two natural
transformation the tensor product $\eta\ot \phi:F_{X\ot Y}\to
F_{Z\ot W}$ is given by the composition
\begin{align}\label{prod-morf}(\eta\ot \phi)_M= m^{-1}_{ZWM}\eta_{W\otb M} (\id_X\ot\phi_M)
m_{XYM},
\end{align}
for all $M\in \Mo$. The unit element is $F_{\uno}$.

\begin{rmk} Note that for any $X, Y, U,
V\in\ca$ and $f:X\to Y, g:U\to V$,
\begin{align} (\eta_f\ot  \eta_g)_M= (f\ot g)\ot \id_M,
\end{align}
for all $M\in \Mo$.
\end{rmk}

\begin{prop}  If $\Mo\simeq \No$ as module categories then there is
    a tensor equivalence
    $\Mo\ltimes \ca\simeq \No\ltimes \ca$.

\end{prop}
\pf  Assume that $(\Fc, c):\Mo\to \No$ and $(\Gc, d):\No\to \Mo$
is a pair of equivalence of module categories. Let $\theta: \Id
\to \Fc\circ \Gc$ be a natural isomorphism of module functors,
that is $\theta$ satisfies
\begin{align}\label{isonatural} c_{X,\Gc(N)}\, \Fc(d_{X,N})\,
\theta_{X\otb N}=\id_X\ot \theta_N,
\end{align}
for all $X\in \ca$, $N\in \No$.

\medbreak

Define $\Phi:\Mo\ltimes \ca\to \No\ltimes \ca$ the functor
$\Phi(F_X)=\widetilde{F}_X$, for any $X\in \ca$. Here, we denote
$\widetilde{F}_X:\No\to \No$ the functor $\widetilde{F}_X(N)=X\otb
N$, for all $N\in \No$. If $X, Y\in \ca$ and $\eta:F_X \to F_Y$ is
a natural transformation then $\Phi(\eta):\widetilde{F}_X\to
\widetilde{F}_Y$ is given by the composition

\begin{align}\label{comp}\begin{split} &X\otb N\xrightarrow{\id_X\ot \theta_N} X\otb
\Fc(\Gc(N))\xrightarrow{c^{-1}_{X,\Gc(N)}} \Fc(X\otb
\Gc(N))\xrightarrow{}\\
&\xrightarrow{\Fc(\eta_{\Gc(N)})} \Fc(Y\otb \Gc(N))
\xrightarrow{\Fc(d^{-1}_{Y,N})} \Fc\Gc(Y\otb
N)\xrightarrow{\theta^{-1}_{Y\otb N}} Y\otb N,
\end{split}
\end{align}
for all $N\in \No$. The tensor structure on $\Phi$ is given by the
identity. That is, for any $X, Y\in \ca$ the natural isomorphisms
$\xi:\Phi(F_X\ot F_Y)\to \Phi(F_X)\ot \Phi(F_Y)$ are given
$$ \xi_{X,Y,N}:\widetilde{F}_{X\ot Y}(N)\to  \widetilde{F}_{X\ot Y}(N),
\quad \xi_{X,Y,N}=\id_{X\ot Y}\ot \id_N,$$ for all $N\in \No$. We
have to check that for all $X, Y, Z\in \ca$, $N\in\No$ the
following identity holds:
\begin{align}\label{t11} (a_{XYZ}\ot \id_N) (\xi_{XY}\id_Z)_N\; \xi_{X\ot Y,Z, N}= (\id_X\ot\xi_{YZ})_N
\Phi(a_{XYZ}\ot \id_N).\end{align}

The left hand side of \eqref{t11} is equal to  $(a_{XYZ}\ot
\id_N)$, the right hand side is equal to $\Phi(a_{XYZ}\ot \id_N)$,
and using \eqref{comp}, is  equal to
\begin{align*} &\theta^{-1}_{(X\ot(Y\ot Z))\otb N}
\; \Fc(d^{-1}_{X\ot(Y\ot Z), N})\; (a_{XYZ}\ot\id_{\Fc(\Gc(N))})
c^{-1}_{(X\ot
Y)\ot Z, \Gc(N)} (\id\ot\theta_N)\\
&=(a_{XYZ}\ot\id_N) \theta^{-1}_{((X\ot Y)\ot Z)\otb N} \;
\Fc(d^{-1}_{(X\ot Y)\ot Z, N})c^{-1}_{(X\ot Y)\ot Z, \Gc(N)}
(\id\ot\theta_N)\\
&=(a_{XYZ}\ot\id_N).
\end{align*}
The last equality follows by \eqref{isonatural}.

\epf

The following definition seems to be well-known.
\begin{defi} A \emph{cocycle} in $\ca$ is a family of
isomorphisms $J_{X,Y}\in \Hom_{\ca}(X\ot Y)$ such that for all $X,
Y, Z\in \ca$
\begin{align}\label{cocyc1} a_{XYZ}\, J_{X\ot Y,Z}\,(J_{X,Y}\ot\id_Z)=
J_{X,Y\ot Z} \,(\id_X\ot J_{Y,Z})\, a_{XYZ},
\end{align}
\begin{align}\label{cocyc2} J_{X,\uno}=\id_{X\ot \uno}, \quad
J_{\uno,X}=\id_{\uno\ot X }.
\end{align}

\end{defi}
If $J$ is a cocycle in $\ca$ then there is a new monoidal
category, $\ca^J$ defined as follows. The objects and morphisms
are the same as in $\ca$. The tensor product of $\ca^J$ coincides
with the tensor product of $\ca$ on objects. If $f:X\to Y, g:Z\to
W$ is a pair of morphisms then
$$f\widetilde{\ot} g= J_{Y,W} (f\ot g) J^{-1}_{X,Z}. $$
Evidently if $J$ commutes with morphisms in $\ca$ then the tensor
category $\ca^J$ is equivalent to $\ca$.

Let $(\Mo, m, l)$ be a module category over $\ca$. The following
definition is due to Donin and Mudrov, see \cite[Definition
5.2]{DM1}.
\begin{defi} A \emph{dynamical twist} for the extension $\Mo\ltimes
\ca$ is a cocycle $J$ in $\Mo\ltimes \ca$ such that $J$ commutes
with morphisms in $\ca$, that is
\begin{align}\label{dynam1} J_{Z,W} (\eta_f\ot\, \eta_g)= (\eta_f\ot\, \eta_g)
J_{X,Y},
\end{align}
for all $f\in \Hom_{\ca}(X,Z), g\in \Hom_{\ca}(Y,W)$.
\end{defi}

More explicitly, a dynamical twist is a family of isomorphisms
$$ J_{X,Y,M}:(X\ot Y)\otb M\to (X\ot Y)\otb M,$$
$X, Y\in\ca, M\in \Mo$ such that
\begin{align}\label{dynamic1} \begin{split}(a_{XYZ}\ot\id_M)&J_{X\ot Y,Z,M} m^{-1}_{X\ot Y, Z,M}
J_{X,Y,Z\otb M} m_{X\ot Y,Z,M}=\\
&=J_{X,Y\ot Z, M}m^{-1}_{X,  Y\ot Z,M}(\id_X\ot J_{Y,Z,M}) m_{X,
Y\ot Z,M} (a_{XYZ}\ot\id_M),
\end{split}
\end{align}
\begin{align}\label{dynamic2} (l_X\ot \id_M)J_{X, \uno, M}=(l_X\ot
\id_M), \;\; J_{\uno, X , M}(r_X\ot \id_M)=(r_X\ot \id_M),
\end{align}
for all $X, Y, Z\in \ca$ , $M\in \Mo$.

\medbreak

Equation \eqref{dynam1} implies that
$$J_{Z,W, M}m^{-1}_{Z,W,M}(f\ot (g\ot\id_M))m_{X,Y,M}=
m^{-1}_{Z,W,M}(f\ot (g\ot \id_M))m_{X,Y,M} J_{X,Y, M},$$ for all
morphisms $f:X\to Z$, $g:Y\to W$ in $\ca$, and all $M\in \Mo$.

\subsection{Module categories coming from dynamical twists}

If $J$ is a dynamical twist for the dynamical extension
$\Mo\ltimes \ca$ we will denote by $\Mo^{(J)}$ the category $\Mo$
with the following module category structure; the action is the
same as in $\Mo$ and the associativity isomorphisms are
$$\widehat{m}_{X,Y,M}:(X\ot Y)\otb M\to X\otb ( Y\otb M), \quad
\widehat{m}_{X,Y,M}= m_{X,Y,M} J^{-1}_{X,Y,M},$$ for all $X, Y\in
\ca$, $M\in \Mo$.

\begin{prop}\label{dyn-modcat} Let $J$ be a dynamical twist for
the dynamical extension $\Mo\ltimes \ca$.\begin{itemize}
    \item[1.]  $\Mo^{(J)}$ is a module category over $\ca$. If
    $\Mo$ is indecomposable then so is $\Mo^{(J)}$.
    \item[2.] There is a tensor equivalence $\Mo^{(J)}\ltimes \ca
    \simeq (\Mo\ltimes \ca)^J$.
\end{itemize}
\end{prop}
\pf Straightforward. \epf

The idea of using the module category language in the study of
dynamical twists is due to Ostrik, see \cite{O1}. In \emph{loc.
cit.} the author relates the classification of module categories
over $\Rep(G)$, $G$ a finite group, with the results obtained by
Etingof and Nikshych on dynamical twists over the group algebra of
the group $G$ \cite{EN1}. This idea was used with profit in
\cite{M}.

\subsection{ Dynamical twists and dynamical adjoint functors}


\medbreak

Let $\ca$ be a tensor category. Cocycles in $\ca$ are in bijective
correspondence with natural associative operations in the space
$\Hom_{\ca}$. If $J$ is a cocycle then $$ \circledast:
\Hom_{\ca}(V,U)\ot \Hom_{\ca}(V',U')\to \Hom_{\ca}(V\ot V',U\ot
U')$$ defined by $\phi\circledast\psi= (\phi\ot \psi)J^{-1}_{VU} $
is an associative operation for all $V,U, V',U'\in\ca$.

\medbreak

The following result is analogous to \cite[Lemma 6.1]{DM1}.

\begin{lema}\label{asso} Assume that there is an associative operation
$$\circledast:
\Hom_{\ca}(V,U)\ot_{\ku} \Hom_{\ca}(V',U')\to \Hom_{\ca}(V\ot
V',U\ot U')$$  for all $V,U, V',U'\in\ca$ such that
\begin{align}\label{conm1} (\phi \circledast \psi
)=(\phi\ot \psi)(\id\circledast \id),
\end{align}
\begin{align}\label{conm2} \phi \circledast\xi=\phi\ot \xi, \; \xi\circledast\phi
=\xi\ot\phi
\end{align}
for all morphisms $\phi,\psi,\alpha,\beta$ in $\ca$ and
$\xi\in\Hom_{\ca}(V , \uno)$. Assume also that for any $U,V\in\ca$
$I_{UV}=\id_U\circledast\id_V$ is invertible. Then
$J_{UV}=I^{-1}_{UV}$, $U,V\in\ca$, is a cocycle in $\ca$.

\end{lema}

\pf The proof is entirely similar to the proof of \cite[Lemma
6.1]{DM1}.
 \epf

\medbreak

 Let $\ca, \ca'$ be two tensor categories. Let $(\Mo, m, l)$ be a
 module category over $\ca$ and $(\Mo', m', l')$ be a module category over $\ca'$. Let
$(R,b):\ca\to \ca'$ be a tensor functor. The following definition
is \cite[Def. 6.2]{DM1} for right module categories.

\begin{defi} A functor $T:\Mo'\to \Mo$ is said to be a \emph{dynamical
adjoint} to $R$ if there exists a family of natural isomorphisms
\begin{align*}\xi_{X,M,N}:\Hom_{\Mo}(X\otb\; T(M),
T(N))\xrightarrow{\;\;\simeq\;\;} \Hom_{\Mo'}(R(X)\otb\, M,N),
\end{align*}
\end{defi}
for all $X\in \ca$, $M, N\in \Mo'$. We further assume that for any
$M\in \Mo'$
\begin{align}\label{req1} \xi_{\uno, M,M}(l_{T(M)})=l'_M.
\end{align}

\begin{rmk}
For any $M, M', N, N'\in \Mo'$ and $X, Y\in \ca$ and morphisms
$f:N\to N'$, $g:M:\to M'$, $\alpha:X\to Y$, the naturality of
$\xi$ implies that

\begin{align}\label{natxi1} f\circ \xi_{X,M,N}(\beta_1)=
\xi_{X,M,N'}(T(f)\beta_1),
\end{align}
\begin{align}\label{natxi2} \xi_{X,M',N}(\beta_2)\circ(\id_{R(X)}\ot
g)=\xi_{X,M,N}(\beta_2(\id_X\ot T(g))),
\end{align}
\begin{align}\label{natxi3} \xi_{Y,M,N}(\beta_3)(R(\alpha)\ot
\id_M)= \xi_{X,M,N}(\beta_3(\alpha\ot\id_{T(M)})),
\end{align}
for any $\beta_1\in \Hom_{\Mo}(X\otb \,T(M), T(N)),\,$ $\beta_2\in
\Hom_{\Mo}(X\otb \,T(M'), T(N))\,$, $\beta_3 \in \Hom_{\Mo}(Y\otb
T(M), T(N)).$

\end{rmk}

\medbreak

The category $\Mo'$ is a module category over $\ca$ via $R$. The
action is given by $\otc:\ca\times \Mo'\to \Mo'$, $X\otc
M=R(X)\otb M$, for all $X\in \ca$, $M\in\Mo'$. The associativity
isomorphisms are
$$ \widetilde{m}_{X,Y,M}= m'_{R(X),R(Y),M} (b_{XY}\ot\id_M)
$$
for all $X, Y\in \ca$, $M\in\Mo'$.

\bigbreak

For any pair of dynamical adjoint functors $(R,T)$ we will repeat
the construction given in \cite{DM1} of a dynamical twist for the
extension $\Mo'\ltimes \ca$. For this we will define an
associative operation in $\Hom_{\Mo'\ltimes \ca}$.

\smallbreak

In some sense, the dynamical twist constructed from the pair
$(R,T)$ \emph{measures} how far is the functor $T$ from being a
module functor.

\medbreak

Let $X, Y, U,V\in \ca$ and $\phi:F_X\to F_Y$, $\psi:F_U\to F_V$ be
morphisms in $\Mo'\ltimes \ca$. So for each $M\in \Mo'$ we have
that $\phi_M: R(X)\otb M\to R(Y)\otb M$, $\psi_M:R(U)\otb M\to
R(V)\otb M$ are morphisms in $\Mo'$.

\medbreak

Set $\widetilde{\phi}_M=\xi^{-1}_{X,M,R(Y)\otb M}(\phi_M)$, and
$\widetilde{\psi}_M=\xi^{-1}_{U,M,R(V)\otb M}(\psi_M)$. Thus we
define $(\phi\circledast \psi)_M$ as the image by $\xi$ of the
composition
$$ (X\ot U)\otb T(M)\xrightarrow{\;\;m_{X,U,T(M)}\;\;}X\otb(U\otb T(M))
 \xrightarrow{\;\id_X\ot \widetilde{\psi}_M\;} X\otb
 T(R(V)\otb M)\to
$$
$$\xrightarrow{\;\;\widetilde{\phi}_{R(V)\otb M}\;\;}
T(R(Y)\otb (R(V)\otb
M))\xrightarrow{\;\;T(\widetilde{m}^{-1}_{YVM})\;\;}T(R(Y\ot
V)\otb M).$$

That is,  for all $M\in \Mo'$, $(\phi\circledast \psi)_M$ equals
to
$$\xi_{X\ot U, M,R(Y\ot V)\otb M}\big(T(\widetilde{m}^{-1}_{YVM})
\widetilde{\phi}_{R(V)\otb M}(\id_X\ot \widetilde{\psi}_M)
m_{X,U,T(M)}\big).$$

\begin{lema}\label{conmutacion}  For any $X, Y, U, V\in \ca$ and
 morphisms $f:X\to Y$, $g\in U\to V$, $\phi:F_X\to F_Y$,
 $\psi:F_U\to F_V$
\begin{align}\label{conmutacion1} (\eta_f \circledast\eta_g)=
(\id_Y\circledast\id_V) (\eta_f\ot \eta_g).
\end{align}
\begin{align}\label{conmutacion2} \phi\circledast\psi=
(\phi\ot \psi)(\id_X\circledast\id_U).
\end{align}
\end{lema}

\pf Using \eqref{natxi3} we have that $(\eta_f \circledast\eta_g)$
equals to
\begin{align*}\hspace{-0.9cm}=\xi_{X\ot U, M,R(Y\ot V)\otb
M}&\big(T(\widetilde{m}^{-1}_{YVM})\xi^{-1}_{Y,M,R(Y)\otb M}(\id)
(f\ot\id_{T(M)})\\ & (\id_X\ot \xi^{-1}_{V,M,R(V)\otb
M}(\id)(g\ot\id_{T(M)})) m_{X,U,T(M)}\big)
\end{align*}
\begin{align*}=\xi_{X\ot U, M,R(Y\ot V)\otb
M}&\big(T(\widetilde{m}^{-1}_{YVM})\xi^{-1}_{Y,M,R(Y)\otb M}(\id)
(\id_Y\ot \xi^{-1}_{V,M,R(V)\otb M}(\id))\\ &(f\ot
(g\ot\id_{T(M)})) m_{X,U,T(M)}\big)
\end{align*}
\begin{align*}=\xi_{X\ot U, M,R(Y\ot V)\otb
M}&\big(T(\widetilde{m}^{-1}_{YVM})\xi^{-1}_{Y,M,R(Y)\otb M}(\id)
(\id_Y\ot \xi^{-1}_{V,M,R(V)\otb M}(\id))\\ &m_{Y,V,T(M)}((f\ot
g)\ot\id_{T(M)}) \big)
\end{align*}
\begin{align*}=\xi_{X\ot U, M,R(Y\ot V)\otb
M}&\big(T(\widetilde{m}^{-1}_{YVM})\xi^{-1}_{Y,M,R(Y)\otb M}(\id)
(\id_Y\ot \xi^{-1}_{V,M,R(V)\otb M}(\id))\\ &m_{Y,V,T(M)}
\big)(R(f\ot g)\ot\id_{M})=(\id_Y\circledast\id_V) (\eta_f\ot
\eta_g).
\end{align*}
The third equality follows from the naturality of $m$ and the
fourth equality, again, follows from \eqref{natxi3}. Thus we have
proved \eqref{conmutacion1}.

\smallbreak

Note that for any $\phi:F_X\to F_Y$, $\psi:F_U\to F_V$, equations
\eqref{natxi1} and \eqref{natxi2} implies that for any $M\in\Mo'$
$\widetilde{\phi}_{R(V)\otb M}(\id_X\ot \widetilde{\psi}_M)$
equals to
$$ \xi^{-1}_{X, R(U)\otb M,R(Y)\otb(R(V)\otb M) }\big(
\phi_{R(V)\otb M}(\id_X\ot \psi_M)\big)\, \big(\id_X\ot \,
\xi^{-1}_{U,M,R(U)\otb M}(\id) \big). $$ Also, using
\eqref{natxi1} we get that $\xi^{-1}_{X, R(U)\otb
M,R(Y)\otb(R(V)\otb M) }\big( \phi_{R(V)\otb M}(\id_X\ot
\psi_M)\big)$ is equal to

$$T(\phi_{R(V)\otb M}((\id_X\ot \psi_M))\, \xi^{-1}_{X, R(U)\otb M,
R(X)\otb(R(U)\otb M) }(\id).$$

Thus, for any $M\in\Mo'$, $\xi^{-1}_{X\ot U, M,R(Y\ot V)\otb
M}\big((\phi\circledast\psi)_M\big)$ is equal to
\begin{align*}T(\widetilde{m}^{-1}_{Y,V,M}\phi_{R(V)\otb M}((\id_X\ot
\psi_M)) &\xi^{-1}_{X, R(U)\otb M, R(X)\otb(R(U)\otb M) }(\id)
\big (\id_X\ot \\
 &\xi^{-1}_{U,M,R(U)\otb M}(\id) \big)\end{align*}

Follows from \eqref{prod-morf} that
$$ T(\widetilde{m}^{-1}_{Y,V,M}\phi_{R(V)\otb M}((\id_X\ot
\psi_M))= T((\phi\ot \psi)_M\, \widetilde{m}^{-1}_{X,U,M}),$$
hence $(\phi\circledast\psi)_M$ equals to

\begin{align*}\xi_{X\ot U, M,R(Y\ot V)\otb
M}\big( T((\phi\ot \psi)_M\, \widetilde{m}^{-1}_{X,U,M})
&\xi^{-1}_{X, R(U)\otb M, R(X)\otb(R(U)\otb M) }(\id) \big
(\id_X\ot \\
&\xi^{-1}_{U,M,R(U)\otb M}(\id) \big)  \big)
\end{align*}
Using again \eqref{natxi1} we get that $(\phi\circledast\psi)_M$
equals to

\begin{align*}(\phi\ot \psi)_M \, \xi_{X\ot U, M,R(Y\ot V)\otb
M}\big( T( \widetilde{m}^{-1}_{X,U,M}) &\xi^{-1}_{X, R(U)\otb M,
R(X)\otb(R(U)\otb M) }(\id) \big
(\id_X\ot \\
&\xi^{-1}_{U,M,R(U)\otb M}(\id) \big)  \big),
\end{align*}
and by definition this is equal to $(\phi\ot
\psi)_M(\id_X\circledast\id_U)_M.$

 \epf

\begin{defi} For any $X, Y\in \ca$ set
\begin{equation}\label{inv-twist} I_{X,Y}=\id_X \circledast\id_Y.
\end{equation}
\end{defi}

\begin{lema}\label{twist-uni} For any $X\in \ca$, $I_{X,\uno}=\id_X\ot\id_{\uno}$,
$I_{\uno,X}=\id_{\uno}\ot\id_X$.
\end{lema}
\pf By definition, for all $M\in \Mo'$, $\xi^{-1}_{X\ot \uno, M,
R(X\ot\uno)\otb M}((\id_X \circledast \id_{\uno})_M)$ equals to
\begin{align}\label{step1} T(\widetilde{m}^{-1}_{X,\uno,M})
\xi^{-1}_{X, R(\uno)\otb M, R(X)\otb(\uno\otb M)}(\id)(\id_X\ot
\xi^{-1}_{\uno, M,\uno\otb M}(\id))m_{X, \uno, M}
\end{align}
Using \eqref{modcat2} and \eqref{natxi1} for the map
$\id_{R(X)}\ot l_M$ we get that \eqref{step1} is equal to
\begin{align}\label{step2} T(b_{X,\uno}r'^{-1}_{R(X)}\ot\id_M)\xi^{-1}_{X, \uno\otb
M, R(X)\otb M}(\id_{R(X)}\ot l'_M) (\id_X\ot \xi^{-1}_{\uno,
M,\uno\otb M}(\id))m_{X, \uno, M}.
\end{align}

From \eqref{natxi1} and \eqref{req1} we get that
\begin{align*}\xi^{-1}_{\uno, M,\uno\otb
M}(\id)&=T(l'^{-1}_M)\xi^{-1}_{\uno,M,M}(l'_M)=T(l'^{-1}_M)l_{T(M)}.
\end{align*}
Now, using \eqref{modcat2} follows that \eqref{step2} equals to

\begin{align}\label{step3} T(b_{X,\uno}r'^{-1}_{R(X)}\ot\id_M)\xi^{-1}_{X, \uno\otb
M, R(X)\otb M}(\id_{R(X)}\ot l'_M)(\id_X\ot
T(l'^{-1}_M))(r_X\ot\id_M).
\end{align}

From \eqref{natxi2} and  \eqref{natxi2} we get that
$$\xi^{-1}_{X, \uno\otb
M, R(X)\otb M}(\id_{R(X)}\ot l'_M)(\id_X\ot
T(l'^{-1}_M))=\xi^{-1}_{X,M,R(X)\otb M}(\id), $$
$$\xi^{-1}_{X,M,R(X)\otb M}(\id)
(r_X\ot\id_M)= \xi^{-1}_{X\ot\uno, M,R(X)\otb M}(R(r_X)\ot\id_M)$$
thus \eqref{step3} is equal to
\begin{align*}\label{step4}  T(b_{X,\uno}r'^{-1}_{R(X)}\ot\id_M)&\xi^{-1}_{X,M,R(X)\otb M}(\id)
(r_X\ot\id_M)=\\
&=T(b_{X,\uno}r'^{-1}_{R(X)}\ot\id_M)\xi^{-1}_{X\ot\uno,
M,R(X)\otb M}(R(r_X)\ot\id_M).
\end{align*}
Finally, $(\id_X \circledast \id_{\uno})_M$ is equal to
\begin{align*}&=\xi_{X\ot \uno, M,
R(X\ot\uno)\otb
M}\big(T(b_{X,\uno}r'^{-1}_{R(X)}\ot\id_M)\xi^{-1}_{X\ot\uno,
M,R(X)\otb M}(R(r_X)\ot\id_M)\big)\\
&=(b_{X,\uno}r'^{-1}_{R(X)}\ot\id_M)\xi_{X\ot\uno, M,R(X)\otb
M}\big(\xi^{-1}_{X\ot\uno, M,R(X)\otb M}(R(r_X)\ot\id_M)\big)\\
&=(b_{X,\uno}r'^{-1}_{R(X)}\ot\id_M)(R(r_X)\ot\id_M)=\id_{R(X)}\ot\id_{\uno}\ot\id_M.
\end{align*}
The equality $I_{\uno,X}=\id_{\uno}\ot\id_X$ follows in a similar
way. \epf

\begin{lema}\label{asociat} For any $U, V, W, X, Y, Z\in \ca$, $\phi:F_X\to F_Y$,
$\psi:F_U\to F_V$, $\chi:F_Z\to F_W$
$$(\phi\circledast\psi)\circledast\chi=\phi\circledast(\psi \circledast\chi).$$
\end{lema}

\pf The proof is done by a tedious but straightforward
computation.

\epf

The following result is a re-statement of \cite[Prop. 6.3]{DM1}.

\begin{teo}\label{dyt} Let us assume that for any $X\in \ca, M\in \Mo'$ the
map $\xi^{-1}_{X, M,R(X)\otb M}(\id)$ is an isomorphism. Then for
any $X,Y\in \ca$ the maps $I_{X,Y}$ are invertible and
$J_{X,Y}=I^{-1}_{X,Y}$ is a dynamical twist for the extension
$\Mo'\ltimes \ca$.
\end{teo}

\pf Once we prove that $I_{X,Y}$ are invertible for any $X,Y\in
\ca$, the proof that $I^{-1}_{X,Y}$ is a dynamical twist for the
extension $\Mo'\ltimes \ca$ follows immediately from Lemmas
\ref{asso}, \ref{conmutacion}, \ref{twist-uni} and \ref{asociat}.

\medbreak

Let $X\in \ca$, $M, N\in \Mo'$ and $\beta:X\ot T(M)\to T(N)$ be an
invertible map with inverse $\gamma:T(N)\to X\ot T(M)$. There
exists an $f: R(X)\otb M\to N$ such that
$\beta=\xi^{-1}_{X,M,N}(f)$. Using \eqref{natxi1} we obtain that
$\beta=T(f)\xi^{-1}_{X, M,R(X)\otb M}(\id)$, hence $f$ is
invertible. Therefore $\xi_{X, M,N}\big(\beta\big)$ is invertible.
In another words, $\xi$ maps isomorphisms to isomorphisms.

\medbreak

Let $X,Y\in \ca$, $M\in \Mo'$. By definition
$$ (I_{X,Y})_M=\xi_{X\ot Y, M, R(X\ot Y)\otb
M}\big(T(\widetilde{m}^{-1}_{X,Y,M}) \;\theta\;
m_{X,Y,T(M)}\big),$$ where $\theta=\xi^{-1}_{X, R(Y)\otb M,
R(X)\otb (R(Y)\otb M)}(\id)\, (\id_X\ot \xi^{-1}_{Y,M,R(Y)\otb
M}(\id)).$ By assumption $\theta$ is invertible, thus
$T(\widetilde{m}^{-1}_{X,Y,M}) \;\theta\; m_{X,Y,T(M)}$ is
invertible, hence $ I_{X,Y}$ is invertible.
 \epf

\begin{rmk} For any $X,Y\in \ca$  $ I_{X,Y}=\id_{X\ot Y}$ if and
only if the functor $(T, c):\Mo'\to \Mo$ is a module functor,
where $c_{X,M}:T(R(X)\otb M)\to X\otb T(M)$ is defined by
$c_{X,M}=(\xi^{-1}_{X, M,R(X)\otb M}(\id))^{-1}$ for any $X\in
\ca$, $M\in \Mo'$.
\end{rmk}

\begin{lema} Let $T$ be a dynamical adjoint to $R$ and let $J$ be
the dynamical twist associated. Then the functor $(T,
c):(\Mo')^{(J)}\to \Mo$ is a module functor, where
$c_{X,M}=(\xi^{-1}_{X, M,R(X)\otb M}(\id))^{-1}$ for any $X\in
\ca$, $M\in \Mo'$.\qed
\end{lema}

\section{Dynamical twists over Hopf algebras}

In this section we shall focus our attention to the computation of
dynamical twists for a dynamical extension of the category
$\Rep(H)$, where $H$ is a finite-dimensional Hopf algebras, and we
shall give explicit examples.

\bigbreak

Hereafter we shall denote by $H$ a finite-dimensional Hopf
algebra. Let $S$ be a left $H$-comodule algebra.

\begin{defi}\label{defi-dyn} A \emph{dynamical twist with base} $S$ for $H$ is an
invertible element $J\in H\ot_{\ku} H\ot_{\ku} S$ such that
\begin{equation}\label{dynt0}  J^1 s\_{-2} \ot
 J^2 s\_{-1}\ot J^3 s\_{0} =s\_{-2} J^1 \ot
 s\_{-1}J^2 \ot s\_{0}J^3\quad \text{ for all } s\in S,
\end{equation}
\begin{equation}\label{dynt1} j^1\_1 J^1\ot j^1\_2 J^2\ot j^2
J^3\_{-1}\ot j^3 J^3\_0=j^1\ot j^2\_1J^1\ot j^2\_2 J^2\ot j^3 J^3,
\end{equation}
\begin{equation}\label{dynt2} \langle\varepsilon , J^1\rangle \,J^2\ot J^3= 1_H\ot 1_K=
\langle\varepsilon , J^2\rangle \,J^1\ot J^3.
\end{equation}
\end{defi}
Here we use the notation $J=J^1\ot J^2\ot J^3=j^1\ot j^2\ot j^3$
avoiding the summation symbol.

\medbreak

\begin{rmk} Definition \ref{defi-dyn} coincides with the
definition of dynamical twist over an Abelian group given in
\cite{EN1}, \cite{EN2}.

\end{rmk}

\begin{defi}  Two dynamical twists for $H$ over
$S$, $J$ and $J'$, are said to be \emph{gauge equivalent} if there
exists an invertible element $t\in H\ot S$ such that

\begin{align*}\langle \varepsilon , t^1 \rangle\; t^2&= 1,\\
 t^1\_1 J^1\ot t^1\_2 J^2\ot t^2 J^3 &= j^1t^1\ot  j^2
T^1t^2\_{-1} \ot j^3 T^2t^2\_{0}.
\end{align*}
Here $J=J^1\ot J^2\ot J^3$, $J'=j^1\ot j^2\ot j^3$, $t=t^1\ot
t^2=T^1\ot T^2$.
\end{defi}

The following Lemma is a straightforward consequence of the
definitions.

\begin{lema} Let $J$ be a dynamical twist with base $S$ for $H$.
For any $X, Y\in \Rep(H)$, $M\in  {}_S\Mo$ define
$\mathcal{J}_{X,Y,M}: X\ot_{\ku} Y \ot_{\ku} M\to  X\ot_{\ku} Y
\ot_{\ku} M$ by
$$\mathcal{J}_{X,Y,M}(x\ot y\ot m)= J^{1}\cdot x\ot J^{2}\cdot y\ot
J^{3}\cdot m,$$ for all $x\in x, y\in Y, m\in M$. Then
$\mathcal{J}$ is a dynamical twist for the dynamical extension
${}_S\Mo\ltimes \Rep(H)$.

Moreover, any dynamical twist for the extension ${}_S\Mo\ltimes
\Rep(H)$ comes from a dynamical twist with base $S$ over $H$.
\end{lema}\qed

\begin{lema}\label{mod-y-twist} If $J$ and $J'$ are gauge equivalent dynamical twists
for $H$ with base $S$ then ${}_S\Mo^{(J)}\simeq {}_S\Mo^{(J')}$ as
module categories over $\Rep(H)$.
\end{lema}

\pf Let $(F,c):{}_S\Mo^{(J)}\to {}_S\Mo^{(J')}$ be the functor
defined as follows. For any $M\in {}_S\Mo$ and $X\in \Rep(H)$,
$F(M)=M$ and $c_{X,M}:X\otk M\to  X\otk M$, $c_{X,M}(x\ot m)=
t^1\cdot x\ot\, t^2\cdot m$, for any $x\in X, m\in M$. It is easy
to verify that $(F,c)$ gives an equivalence of module categories.

\epf

\subsection{Dynamical twists and Hopf Galois extensions}
Let $S$ be a left $H$-comodule algebra. For any  dynamical twist
with base $S$ for $H$ there is associated a $H$-Galois extension
with coinvariants $S$.

\medbreak

Set $B=H^*\otk S$. The coproduct of $H^*$ endows $B$ with a right
$H^{*\cop}$-comodule structure, that is $\delta:B\to B\otk H^*$,
$\delta(\alpha\ot s)= \alpha\_2\ot s\ot\alpha\_1$, for all
$\alpha\in H^*$, $s\in S$. Clearly $B^{\co H^{*\cop}}=S$.

\smallbreak

If $J\in H\otk H\otk S$ we endowed $B$ with the following product:

\begin{align}\label{productj} (\alpha\ot k)(\beta\ot s)=(J^1\rightharpoonup
\alpha)(J^2k\_{-1}\rightharpoonup \beta )\ot J^3 k\_0 s,
\end{align}
for all $\alpha, \beta\in H^*$, $k, s\in S$.

\begin{prop}\label{twist-galois} Assume that $J\in H\otk H\otk S$
is a dynamical twist with base $S$, then $B\supset S$ is a
$H^{*\cop}$-Hopf Galois extension.
\end{prop}

\pf First we prove that $B$ is an associative algebra with the
product described in \eqref{productj}. Let $\alpha, \beta,
\gamma\in H^*$, $k, s,r\in S$, then
\begin{align*} &\big((\alpha\ot k)(\beta\ot s)\big)(\gamma\ot r)=\big((J^1\rightharpoonup
\alpha)(J^2k\_{-1}\rightharpoonup \beta )\ot J^3 k\_0
s\big)(\gamma\ot r)\\
&= \!(j^1\_1J^1\!\rightharpoonup\!
\alpha)(j^1\_2J^2k\_{-2}\!\rightharpoonup \beta )
(j^2J^3\_{-1}k\_{-1} s\_{-1}\rightharpoonup \!\gamma)\ot
j^3J^3\_0k\_0s\_0 r.
\end{align*}
On the other hand
\begin{align*} &(\alpha\ot k)\big((\beta\ot s)(\gamma\ot r)\big)=
(\alpha\ot k)\big((J^1\rightharpoonup \beta)
(J^2s\_{-1}\rightharpoonup \gamma)\big)
\ot J^3s\_0r\\
&=(j^1\rightharpoonup \alpha)(j^2\_1k\_{-1}J^1\rightharpoonup
\beta ) (j^2\_2 k\_{-2}J^2 s\_{-1}\rightharpoonup
\gamma)\ot j^3k\_0 J^3s\_0r\\
&=(j^1\rightharpoonup \alpha)(j^2\_1J^1k\_{-1}\rightharpoonup
\beta ) (j^2\_2 J^2k\_{-2} s\_{-1}\rightharpoonup \gamma)\ot
j^3J^3k\_0 s\_0r.
\end{align*}
The last equality follows by \eqref{dynt0}. From \eqref{dynt1}
follows that
$$\big((\alpha\ot k)(\beta\ot s)\big)(\gamma\ot r)=(\alpha\ot k)\big((\beta\ot s)(\gamma\ot
r)\big).$$ The proof that $B$ is a $H^{*\cop}$-comodule algebra is
straightforward. Let $can:B\ot_S B\to B\otk H^{*\cop}$ be the
canonical map; that is $can(a\ot b)= ab\_0\ot b\_1$, for all
$a,b\in B$. In this case,
$$ can((\alpha\ot k)\ot (\beta\ot s))= (J^1\rightharpoonup
\alpha)(J^2k\_{-1}\rightharpoonup \beta\_2 )\ot J^3 k\_0 s \ot
\beta\_1,$$ for all  $\alpha, \beta\in H^*$, $k, s\in S$. It is
easy to see that the inverse of $can$ is given by $can^{-1}:B\otk
H^{*\cop}\to B\ot_S B$,
$$can^{-1}(\gamma\ot r\ot \beta)=J^{-1} \rightharpoonup(\gamma \Ss(\beta\_2))\ot 1\ot
(J^{-2} \rightharpoonup\beta\_1)\ot J^{-3} r,$$ for all $\gamma,
\beta \in H^*$, $r\in S$. Thus $B\supset S$ is a
$H^{*\cop}$-Galois extension.

\epf
\subsection{Dynamical twists coming from dynamical datum}

\

In this subsection we shall give a method for constructing
dynamical twists with base $A$, where $A\subset H$ is a Hopf
subalgebra. This construction is based on the same ideas contained
in \cite{EN1}, \cite{M} without assuming commutativity nor
cocommutativity of the base of the dynamical twist.

\medbreak The following definition generalizes \cite[Def.
4.5]{EN1}, see also \cite[Def. 3.8]{M}.

\begin{defi} A \emph{dynamical datum} for $H$ over $A$ is a pair $(K, T)$ where $K$
is a left $H$-comodule algebra $H$-simple, with trivial
coinvariants, $T:\Rep(A)\to {}_K\Mo$ is a functor such that there
are natural $H$-module isomorphisms

\begin{align} \omega_{VW}:\St_K(T(V),T(W))\xrightarrow{\;\;\simeq\;\;} \big(\Ind^H_A\,
(V\ot_{\ku} W^*)\big)^*,
\end{align}
for any $V, W\in \Rep(A)$. We shall further assume that
\begin{align}\label{req12} \omega_{VV}(1)(\overline{h\ot v\ot f})=\langle\varepsilon,
h\rangle\langle f , v\rangle,
\end{align} for all $h\in H$, $v\in V, f\in
V^*$. \smallbreak

Two dynamical data $(K,T)$ and $(S,U)$ are \emph{equivalent} if
and only if there exists an element $P\in \, ^{H}\Mo_K$ such that
$S\simeq \End_K(P_K)$ as $H$-comodule algebras, and there exists a
family of natural $K$-module isomorphisms
$$ \phi_V: P\ot_K T(V)\xrightarrow{\;\; \simeq \;\;} U(V),$$
for all $V\in \Rep(A)$.
\end{defi}

\begin{rmk} If $(K,T)$ is a dynamical datum, then for any $V\in
\Rep(A)$
\begin{align}\label{dim-dyn-datum}  \dim A\,  (\dim T(V))^2= \dim K \, (\dim V)^2.
\end{align}
These formula follows straightforward from the definition of
dynamical datum and formula \eqref{dim-stab}.
\end{rmk}

Denote by $R:\Rep(H)\to \Rep(A)$ the restriction functor.

\begin{prop}\label{main1} If $(K, T)$ is a dynamical datum then $T$ is a
dynamical adjoint to $R$.

\end{prop}

\pf Category ${}_K\Mo$ is a module category over $\Rep(H)$ as
explained in section \ref{cathopf}. The category $\Rep(A)$ is a
module category over itself. Let $V, W\in {}_K\Mo$ and $X\in
\Rep(H)$. Then
\begin{align*}&\Hom_{K}(X\ot_{\ku}\; T(V), T(W))\simeq \Hom_{H}(X, \St_K(T(V), T(W)))
\simeq\\& \Hom_{H}(X,\big(\Ind^H_A\, (V\ot_{\ku} W^*))\simeq
\Hom_{H}(\Ind^H_A\, (V\ot_{\ku} W^*), X^* )\simeq\\
&\simeq \Hom_{A}(V\ot_{\ku} W^*, R(X^*))\simeq \Hom_{A}(R(X),
W\ot_{\ku} V^*)\simeq\\
&\simeq \Hom_{A}(R(X)\ot_{\ku} V , W ).
\end{align*}

The first isomorphism follows from Proposition
\ref{stab-properties} (3), the second because  $(K, T)$ is a
dynamical datum and the fourth isomorphism is Frobenius
reciprocity.

Let us denote by $\xi:\Hom_{K}(X\ot_{\ku}\; T(V), T(W))\to
\Hom_{A}(R(X)\ot_{\ku} V , W )$ the composition of the above
isomorphisms. Is clear that $\xi$ satisfies \eqref{req1} since we
requested that the isomorphisms $\omega_{VW}$ satisfy
\eqref{req12}.

\epf

\begin{defi} For any dynamical datum $(K,T)$ we shall denote by
$J_T$ the associated dynamical twists for the Hopf algebra $H$
with base $A$ according to Theorem \ref{dyt}.
\end{defi}

In the following we shall prove that the construction of the
dynamical twist does not depend on the equivalence class of the
dynamical datum.

\begin{prop}\label{invariance1}  Let $(K,T)$ and $(S,U)$ be two
equivalent dynamical data over $A$. Then $J_T$ is gauge equivalent
to $J_S$.
\end{prop}

For the proof we will need first some technical results. From the
hypothesis we know that there exists $P\in {}^{H\!}\Mo_{K}$ such
that $S\simeq \End_K(P_K)$ as $H$-module algebras and  natural
isomorphisms
$$\phi_V: P\ot_K T(V)\xrightarrow{\;\;\; \simeq \;\;\;} U(V),$$
for all $V\in \Rep(A).$ For any For any $X\in \Rep(H)$, $V\in
\Rep(A)$ let us denote by
$$ \xi:\Hom_{K}(X\ot_{\ku}\; T(V), T(W))\to
\Hom_{A}(R(X)\ot_{\ku} V , W ),$$
$$ \zeta:\Hom_{S}(X\ot_{\ku}\; U(V), U(W))\to
\Hom_{A}(R(X)\ot_{\ku} V , W ),$$ the family of natural
isomorphisms constructed in the proof of Proposition \ref{main1}.

\smallbreak

For any $X\in \Rep(H)$, $V\in  \Rep(A)$, $M\in {}_K\Mo$ let us
define
$$\theta_{X,M}:X\otk (P\ot_K M)\to P\ot_K (X\otk M)$$   as
follows. For any $x\in X, p\in P, m\in M$
$$\theta_{X,M}(x\ot (p\ot m))=p\_0\ot\,
\Ss^{-1}(p\_{-1})\cdot x\ot\, m .$$ Let us also define
$$\sigma_{X,V}:X\otk U(V)\to U( X\otk V) $$
as the composition
$$\sigma_{X,V}=\phi_{X\otk V} \big(\id_P\ot\, \xi^{-1}_{X,V,X\otk V}(\id)\big)
\,\theta_{X,T(V)}\,(\id_X\ot\phi^{-1}_{V}). $$ Clearly,
$\sigma_{X,V}$ and $\theta_{X,M}$ are isomorphisms.

\begin{lema}\label{tec4} For any $X, Y\in \Rep(H)$, $V, W\in
\Rep(A),$ $M, N \in {}_K\Mo$ and any morphisms $f:X\to Y$,
$\beta:V\to W$, $g:M\to N$
 \begin{align}\label{theta1} \theta_{X\ot_{\ku} Y,M}&=\theta_{X,Y\ot_{\ku}
M}(\id_X\ot\, \theta_{X,M}),\\
\label{theta2} (\id_P\ot\id_X\ot g)\, \theta_{X, M}&=
\theta_{X,N}\, (\id_X\ot\id_P\ot g),\\
\label{theta3} (\id_P\ot f\ot \id_{M}) \, \theta_{X,M}&=
\theta_{Y,M} (f\ot \id_P\ot \id_{M}),
\end{align}
     \begin{align}\label{sigma1} \sigma_{Y,V}(f\ot\id_V)&= U(f\ot\id_V)
    \sigma_{X,V},\\
\label{sigma2} \sigma_{X,W}(\id_X\ot \beta)&= U(\id_X\ot \beta)
\;\sigma_{X,V}.
\end{align}

\end{lema}

\pf Equations \eqref{theta1}, \eqref{theta2} and \eqref{theta3}
are straightforward. By definition
\begin{align*} &\sigma_{Y,V}(f\ot\id_{U(V)})=\phi_{Y\otk V} \big(\id_P\ot\, \xi^{-1}_{Y,V,Y\otk V}(\id)\big)
\,\theta_{Y,T(V)}\,(\id_Y\ot\phi^{-1}_{ V}) (f\ot\id_V)\\
&=\phi_{Y\otk V} \big(\id_P\ot\, \xi^{-1}_{Y,V,Y\otk V}(\id)\big)
\,(\id_P\ot
f\ot\id_{T(V)})\theta_{X,T(V)}\,(\id_Y\ot\phi^{-1}_{V})\\
&=\phi_{Y\otk V} \big(\id_P\ot\, \xi^{-1}_{Y,V,Y\otk
V}(f\ot\id_V)\big) \,\theta_{X,T(V)}\,(\id_Y\ot\,\phi^{-1}_{
V})\\
&=\phi_{Y\otk V} \big(\id_P\ot\, T(f\ot\id_V) \xi^{-1}_{Y,V,Y\otk
V}(\id)\big) \,\theta_{X,T(V)}\,(\id_Y\ot\,\phi^{-1}_{
V})\\
&=U(\id_X\ot \beta)\; \sigma_{X,V}.
\end{align*}
The second equality follows by the naturality of $\phi$ and
\eqref{theta3}, the third equality follows by \eqref{natxi3}, the
fourth by \eqref{natxi1} and the fifth again by the naturality of
$\phi$. Equation \ref{sigma2} follows in a similar way.
\epf

For any $X\in \Rep(H)$, $V\in  \Rep(A)$ set $t_{X,V}:X\otk V\to
X\otk V $ the isomorphism of $A$-modules defined as
$$t_{X,V}=\zeta_{X,V,X\otk V}(\sigma_{X,V}). $$
\begin{lema}\label{tec5} The maps $t_{X,V}$ are natural
isomorphisms. In particular there exists an invertible element
$t=t^1\ot t^2\in H\otk A$ such that for any $x\in X, v\in V$,
$t_{X,V}(x\ot v)= t^1\cdot x\ot\, t^2\cdot v$.
\end{lema}
\pf Let $X, Y\in \Rep(H)$ and let $f:X\to Y$ be any morphism.
\begin{align*} t_{Y,V}\, (f\ot\id_V)&= \zeta_{Y,V,Y\otk V}(\sigma_{Y,V})
\, (f\ot\id_V) =\zeta_{X,V,Y\otk V}
(\sigma_{Y,V}(f\ot\id_{U(V)}))\\&= \zeta_{X,V,Y\otk V}
(U(f\ot\id_V) \,\sigma_{X,V}) =(f\ot\id_V)\;\zeta_{X,V,X\otk V} (
\sigma_{X,V})\\&= (f\ot\id_V)\, t_{X,V}.
\end{align*}
The second equality follows from \eqref{natxi3}, the third by
\eqref{sigma1} and the fourth one by \eqref{natxi1}. The
naturality of $t$ in the second variable follows in an analogous
way using \eqref{sigma2}. \epf

\pf[Proof of Proposition \ref{invariance1}] Let $X, Y\in \Rep(H)$
and $V\in \Rep(A).$ Let $$I_{X,Y,V}, \widetilde{I}_{X,Y,V}:(X\otk
Y)\otk V\to (X\otk Y)\otk V$$ be the isomorphisms defined as
$$I_{X,Y,V}(x\ot y\ot v)= J^{-1}_K\cdot x\ot \,J^{-2}_K\cdot y
\ot \,J^{-3}_K\cdot v, $$
$$ \widetilde{I}_{X,Y,V}(x\ot y\ot v)= J^{-1}_S\cdot x\ot\, J^{-2}_S\cdot y
\ot \, J^{-3}_S\cdot v,$$ for any $x\in X, y\in Y, v\in V$. In
another words, the family of natural isomorphisms $I$ and $
\widetilde{I}$ are given by
$$ I_{X,Y,V}= \xi_{X\otk Y, V, (X\otk Y)\otk V}\big(\xi^{-1}_{X,
Y\otk V, X\otk (Y\otk V)}(\id)\, (\id_X\ot \xi^{-1}_{Y,V,Y\otk
V}(\id)) \big)$$ and
$$ \widetilde{I}_{X,Y,V}= \zeta_{X\otk Y, V, (X\otk Y)\otk V}\big(\zeta^{-1}_{X,
Y\otk V, X\otk (Y\otk V)}(\id)\, (\id_X\ot \zeta^{-1}_{Y,V,Y\otk
V}(\id)) \big).$$

We shall prove that \begin{equation}\label{gauge-equiv}
I_{X,Y,V}\; t_{X\otk Y,V}= t_{X, Y\otk V} (\id_X\ot t_{Y,V}) \,
\widetilde{I}_{X,Y,V}.\end{equation}

Using \ref{natxi1} we obtain that
\begin{align*}\zeta^{-1}_{X\otk Y, V, (X\otk Y)\otk V}\big(I_{X,Y,V}\; t_{X\otk
Y,V}\big)&= U(I_{X,Y,V}) \zeta^{-1}_{X\otk Y, V, (X\otk Y)\otk
V}\big(t_{X\otk Y,V}\big)\\ &= U(I_{X,Y,V})\, \sigma_{X\otk Y,V}.
\end{align*}

Now, using the naturality of $\phi$ and the naturality of $\xi$
\eqref{natxi1} we obtain that $U(I_{X,Y,V})\, \sigma_{X\otk Y,V}$
is equal to
\begin{align*}\phi_{(X\otk Y)\otk V}
\big(\id_P\ot \xi^{-1}_{X\otk Y, V, (X\otk Y)\otk V}
(T(I_{X,Y,V}))\big) \,\theta_{X\otk Y, T(V)} \,(\id_{X\otk
Y}\ot\phi^{-1}_V).
\end{align*}

By definition of the isomorphism $I_{X,Y,V}$ and \eqref{theta1}
this last expression equals to

\begin{align*} \phi_{(X\otk Y)\otk V} \big( \id_P\ot\,\xi^{-1}_{X, Y
\otk V, X\otk (Y\otk V)}(\id)(&\id_X  \ot\, \xi^{-1}_{Y, V, Y\otk
V}(\id)) \big) \theta_{X,Y\ot_{\ku} T(V)}\\&(\id_X\ot\,
\theta_{X,T(V)}) \,(\id_{X\otk Y}\ot\phi^{-1}_V),
\end{align*}
and using \eqref{theta2} equals to

\begin{align*} \phi_{(X\otk Y)\otk V} \big(\id_P &\ot\,\xi^{-1}_{X, Y
\otk V, X\otk (Y\otk V)}(\id) \big) \theta_{X,T(Y\ot_{\ku}
V)}\\&\big(\id_{X\otk P} \ot  \xi^{-1}_{Y, V, Y\otk V}(\id)\big)
(\id_X\ot\, \theta_{X,T(V)}) \,(\id_{X\otk Y}\ot\phi^{-1}_V),
\end{align*}

which is equal to
$$\sigma_{X,Y\otk V}(\id_X\ot\sigma_{Y, V}). $$
From \eqref{natxi1}  follows  that $\zeta^{-1}_{X\otk Y, V, (X\otk
Y)\otk V}\big( t_{X, Y\otk V} (\id_X\ot t_{Y,V}) \,
\widetilde{I}_{X,Y,V}\big) $ \\ equals to
\begin{align*} &=T(t_{X, Y\otk V})\,T(\id_X\ot\,
t_{Y,V})\zeta^{-1}_{X, Y\otk V, X\otk (Y\otk V)}(\id)(\id_X\ot\,
\zeta^{-1}_{Y, V, Y\otk V}(\id))\\
&=T(t_{X, Y\otk V})\,\zeta^{-1}_{X, Y\otk V, X\otk (Y\otk
V)}(\id_X\ot\, t_{Y,V})(\id_X\ot\,
\zeta^{-1}_{Y, V, Y\otk V}(\id))\\
&=T(t_{X, Y\otk V})\zeta^{-1}_{X, Y\otk V, X\otk (Y\otk V)}(\id)
\big(\id_X\ot\, T(t_{Y,V})\zeta^{-1}_{Y, V, Y\otk V}(\id)\big)\\
&=\zeta^{-1}_{X, Y\otk V, X\otk (Y\otk V)}(t_{X, Y\otk V})
\big(\id_X\ot\, \zeta^{-1}_{Y, V, Y\otk V}(t_{Y,V})\big)\\
&=\sigma_{X,Y\otk V}(\id_X\ot\sigma_{Y, V}).
\end{align*}
The second equality follows from \eqref{natxi1}, the third by
\eqref{natxi2} and the fourth again from \eqref{natxi1}. Thus, we
have proven equation \eqref{gauge-equiv}. Follows immediately from
\eqref{gauge-equiv} that the element $t^{-1}$ is a gauge
equivalence for $J_K$ and $J_S$. \epf

There is a reciprocate construction, that is, for any dynamical
twist over $A$ we can associate a dynamical datum as follows.

\medbreak

Let $J$ be a dynamical twist for the dynamical extension
${}_A\Mo\ltimes \Rep(H)$. By Proposition \ref{dyn-modcat} the
category ${}_A\Mo^{(J)}$ is an exact indecomposable module
category over $\Rep(H)$, therefore by Theorem
\ref{cathopf-resultados} there exists an $H$-simple $H$-comodule
algebra $K$ such that $K^{\co H}$ and ${}_A\Mo^{(J)}\simeq
{}_K\Mo$ as module categories over $\Rep(H)$. Let us denote by
$T:{}_A\Mo^{(J)}\to {}_K\Mo$ such equivalence.

\begin{prop}\label{datum-from-twist} The pair $(K, T)$ as above is
a dynamical datum.
\end{prop}
\pf The proof is entirely analogous to the proof of \cite[Prop.
3.18]{M}. For completeness we write the proof. Let $V, W\in
\Rep(A)$, $X\in \Rep(H)$ then 
\begin{align*} &\Hom_H(X,\St_K(T(V),T(W)) )\simeq \Hom_K(X\otk T(V),
T(W))\simeq \\
&\simeq \Hom_K(T(X\otk V), T(W))\simeq \Hom_A(X\otk V, W) \simeq \Hom_A(R(X), W\otk V^*)\\
& \simeq\Hom_A(V\otk W^*, R(X^*))\simeq \Hom_H(\Ind^H_A\,
(V\ot_{\ku} W^*), X^*)\simeq\\
&\simeq \Hom_H(X, \Ind^H_A\, (V\ot_{\ku} W^*)^*).
\end{align*}
The first isomorphism is a consequence of Proposition
\ref{stab-properties} (3), the sixth isomorphism is Frobenius
reciprocity. Thus, the statement follows from Yoneda's Lemma.

\epf

The construction of the dynamical datum from a dynamical twist is
not canonical but it does not depend on the gauge equivalence
class of the dynamical twist.

\begin{prop} Let $J$, $J'$ two gauge equivalent dynamical twists
and let $(K,T)$ and $(S,U)$ the dynamical data associated as in
Proposition \ref{datum-from-twist}. Then $(K,T)$ is equivalent to
$(S,U)$.
\end{prop}
\pf By construction the functors $T:{}_A\Mo^{(J)}\to {}_K\Mo$ and
$U:{}_A\Mo^{(J')}\to {}_S\Mo$ are equivalences of module
categories over $\Rep(H)$. By Lema \ref{mod-y-twist} the
categories ${}_A\Mo^{(J)}$ and ${}_A\Mo^{(J')}$ are equivalent.
Let $G:{}_K\Mo \to {}_A\Mo^{(J)}$ be the inverse of $T$. The
functor $U\circ G: {}_K\Mo \to {}_S\Mo$ is an equivalence of
module categories, thus Proposition \ref{eq1} implies that there
exists an object $P\in {}^{H\!}\Mo_{K}$ such that $S\simeq
\End_K(P_K)$ as $H$-module algebras and  natural isomorphisms
$U(G(M))\simeq P\ot_K M$ for all $M\in {}_K\Mo$. In particular
there are natural isomorphisms
$$U(V)\simeq U(G(T(V)))\simeq P\ot_K T(V), $$
for all $V\in \Rep(A)$.

\epf

\begin{rmk} It would be interesting to know, for a fixed Hopf algebra
$H$, which module categories are equivalent to ${}_A\Mo^{(J)}$ for
some Hopf subalgebra $A$ and a dynamical twist $J$ with base $A$.
\end{rmk}

\medbreak

\subsection{Some examples }

We shall give concrete examples of dynamical datum and explicit
computations of the corresponding dynamical twist.

\subsubsection{$K$ is an $A$-Galois extension}

\

Let us assume that $K$ is an $A$-Galois extension with trivial
coinvariants. Let us denote by $\gamma: A\to K \ot_{\ku} K$ the
map
$$\gamma(a)= can^{-1} (a\ot 1)=a^{[1]} \ot a^{[2]},$$
for all $a\in A$.

Let us assume that $T:\Rep(A)\to {}_K\Mo$ is a functor such that
for any $V, W\in \Rep(A)$ there are $A$-module isomorphisms
$$ T(W)\ot_{\ku} T(V)^*\xrightarrow{\;\;\simeq\;\;} W\ot_{\ku} V^*.$$
The $A$-module structure on $ T(W)\ot_{\ku} T(V)^*$ is given as
follows. If $w\in T(W)$, $f\in T(V)^*$ and $a\in A$, then
$$a\cdot (w\ot f)= a^{[1]}\cdot w \ot f\cdot a^{[2]},$$
where $(f\cdot k)(v)=f(k\cdot v)$, for any $k\in K, v\in V$. This
is a well defined action, see \cite[Lemma 2.21]{AM}.

\begin{lema}\label{example1} Under the above conditions $(K, T)$
is a dynamical datum.
\end{lema}

\pf For any $V, W\in\Rep(A)$ we have that
\begin{align*} \St_K(T(V),T(W))&\simeq \Hom_A(H,
\Hom_{\ku}(T(V),T(W))) \\&\simeq \Hom_A(H, T(W)\ot_{\ku} T(V)^*)
\simeq \Hom_A(H, W\ot_{\ku} V^*)\\&\simeq \big(\Ind^H_A V\ot_{\ku}
W^* \big)^*.
\end{align*}

The first isomorphism follows from \eqref{stab-galois}, and the
last is Lemma \ref{preliminar1}.

 \epf

\subsubsection{Dynamical twists over $\ax(G,\chi,g)$} \

Let $G$ be a finite group, $g\in Z(G)$ and $\chi:G\to \C^{\times}$
a character such that $n=\mid g\mid=\mid\chi(g)\mid$ and
$\chi^n=1$.

\smallbreak

Let us denote by $H=\ax(G,\chi,g)$ the algebra generated by $x, g$
subject to the relations: $x^n=0$, $xh=\chi(h)\, hx$ for all $h\in
G$. The algebra $\ax(G,\chi,g)$ has a Hopf algebra structure as
follows:
$$\Delta(x)=1\ot x+ x\ot g,\; \Delta(f)=f\ot f,\; \varepsilon(x)=0, \;\varepsilon(f)=1, $$
for all $f\in G$.

These Hopf algebras are a special class of monomial Hopf algebras.
See \cite{CYZ}. \medbreak

Let $\lambda \in \C^{\times}$ and let $F\subseteq G$ be a subgroup
such that $g\in F$. In this case
$\ax(F,\chi,g)\subseteq\ax(G,\chi,g)$ is a Hopf subalgebra. Let us
denote by $\ac(F,\lambda)$ the algebra generated by elements $y$,
${e_h:h\in F}$ subject to the relations
$$y^n=\lambda 1, \;\;\; e_h e_f=e_{hf}, \;\;\; ye_h=\chi(h)\; e_h y,$$ for all
$h,f\in F$.

\begin{lema}\label{tec4} Let us denote $\delta:\ac(F,\lambda)\to \ax(G,\chi,g)
\ot \ac(F,\lambda)$, the map given by
$$\delta(y)=g^{-1}\ot y- xg^{-1}\ot 1,\;\; \delta(e_h)=h\ot e_h,$$
for all $h\in F$. Then $\ac(F,\lambda)$ is a left
$\ax(G,\chi,g)$-comodule algebra with trivial coinvariants.
Moreover, $\ac(F,\lambda)$ is a $\ax(F,\chi,g)$-Galois extension.
\end{lema}
\pf Straightforward. \epf

Let $B\subseteq F$ be a subgroup such that $B\bigcap <g>=\{1\}$.
Let $A$ be the group algebra of the group $B$.

Let $T:\Rep(A)\to {}_{\ac(F,\lambda)}\Mo$ be the functor defined
as follows. For any left $A$-module $V$,
$T(V)=\bigoplus_{i=0}^{n-1}\; \C v_i\ot V$. The action of
$\ac(F,\lambda)$ on $T(V)$ is the following. For any $v\in V$,
$i=0\dots n-1$, $h\in F$
$$ y\cdot (v_i\ot v)= \mu\; v_{i-1}\ot v, \;\;\;\; y\cdot (v_0\ot v)= \mu\; v_{n-1}\ot
v$$ $$g \cdot (v_i\ot v)= \chi^i(g)\; v_i\ot v,  \;\;\;\; e_h\cdot
(v_i\ot v)= \chi^i(h)\; v_i\ot h\cdot v.$$

Here $\mu^n=\lambda$ is a fixed $n$-th root of $\lambda$.

\begin{prop}\label{ex2}  The pair $(\ac(F,\lambda), T)$ is a
dynamical datum for $H$ over $A$.
\end{prop}

\pf Let $V, W\in \Rep(A)$. Since $\ac(F,\lambda)$ is
$\ax(F,\chi,g)$-Galois, using \eqref{stab-galois}, we obtain that
$$\St_{\ac(F,\lambda)}(T(V),T(W))\simeq \Hom_{\ax(F,\chi,g)}(\ax(G,\chi,g),
\Hom(T(V), T(W)).$$ Also, by Lemma \ref{preliminar1}, we have that
$$\big( \Ind^{\ax(G,\chi,g)}_{\C B}\, V\ot W^*\big)^*\simeq
\Hom_{\C B}(\ax(G,\chi,g), \Hom(V,W)).$$ Thus, the proof will end
if we prove that $\Hom_{\C B}(\ax(G,\chi,g), \Hom(V,W))$ is
isomorphic to $\Hom_{\ax(F,\chi,g)}(\ax(G,\chi,g), \Hom(T(V),
T(W))$.

\medbreak

Let $G=\bigcup_l F c_l$, $F=\bigcup_{j=0}^{n-1} B g^j$ be right
coset decompositions of $G$ and $F$. Then the algebra
$\ax(G,\chi,g)$ has a basis consisting of elements $\{B g^j
x^ic_l\}_{j,i,l}$.

Let us define the maps $$\phi: \Hom_{\C B}(H,
\Hom(V,W))\xrightarrow{\;\;\;\;} \Hom_{\ax(F,\chi,g)}(H,
 \Hom(T(V), T(W))$$ and $$\psi:\Hom_{\ax(F,\chi,g)}(H,
\Hom(T(V), T(W))\xrightarrow{\;\;\;\;} \Hom_{\C B}(H,\Hom(V,W))$$
defined as follows. If $\xi\in   \Hom_{\C B}(H, \Hom(V,W))$ then
$$\phi(\xi)(c_l)(v_k\ot v)=\sum_{s=0}^{n-1} \, w_s\ot \,\xi(b_sx^kc_l)(v),$$
$$\psi(\alpha)(g^j x^i c_l)(v)=(p_j\ot\id)\big( \alpha(c_l)(v_i\ot v)\big) $$
for any $v\in V$, $k=0\dots n-1$. Here $p_j:\bigoplus_{i=0}^{n-1}
\C w_i\to \C$, $p_j(w_i)=\delta_{ij}$.

It is immediate to verify that these two maps are well defined and
they are one the inverse of each other.

 \epf

\subsection*{Acknowledgments}  This work was partially supported by
 CONICET, Argentina. This work was done during a post-doc fellowship at
Ludwig Maximilians Universit\"{a}t, M\"{u}nchen granted by the
Deutsche Akademische Austauschdienst (DAAD), Germany. The author
thanks L. Feh\'er for providing some references.

\end{document}